\documentclass[11pt, a4paper]{amsart}
\usepackage{amsfonts}
\usepackage{amsmath}
\usepackage{epsfig}
\usepackage{graphicx}
\usepackage{psfrag}
\usepackage{eufrak} 
\usepackage{subfigure}
\newcommand{\equ}[1]{(\ref{#1})} 
 
\newcommand{\be}{\begin{equation}} 
\newcommand{\ee}{\end{equation}} 
\newtheorem{lem}{Lemma}[section] 
\newcommand{\rr}{{\mathbb R}} 
 
\newtheorem{teo}{Theorem}[section] 
\newtheorem{definition}{Definition}[section] 
 
\newtheorem{prop}{Proposition}[section] 
\newtheorem{rem}{Remark}[section] 
 
\newcommand{\calA}{{\mathcal A}} 
 
\newcommand{\calC}{{\mathcal C}} 
\newcommand{\calD}{{\mathcal D}}

\newcommand{\calH}{{\mathcal H}}

\newcommand{\calN}{{\mathcal N}}

\newcommand{\calT}{{\mathcal T}}

\newcommand{\mve}{{n_\epsilon}}

\begin{document} 
\bibliographystyle{plain} 
\title[Uniqueness of self-similar]{Uniqueness of self-similar solutions to the network flow in a given topological class}

\thanks{The author was supported by Proyecto Fondecyt de Iniciaci\'on 11070025}

\author{Mariel S\'aez Trumper}
\address{Mariel S\'aez Trumper
\hfill\break\indent
Edificio Rolando Chuaqui. Facultad de Matem\'aticas
\hfill\break\indent
Pontificia Universidad Cat\'olica de Chile
\hfill\break\indent
Avda. Vicu\~na Mackenna  4860\\
 Macul, Santiago\\
\hfill\break\indent
Chile .}
\email{{\tt  mariel@mat.puc.cl}}

\begin{abstract}
In this paper we study the uniqueness of expanding self-similar solutions to the network flow in a fixed topological class. We prove the result via the parabolic Allen-Cahn approximation proved in \cite{triodginz}.

Moreover, we  prove that any regular evolution of connected tree-like network  (with an initial condition  that might be not regular) 
  is unique in a given a topological class.

\end{abstract}

\maketitle

\section{Introduction}

The flow of planar networks by curve shortening flow has been of interest for several authors in last few years (see \cite{moybycm}, \cite{linsen}, \cite{selsimos} and  \cite{sssnf} for example).  
 A planar network  in $\Omega\subset \rr^2$ is a finite union of arcs embedded on the plane such
 that  each pair of curves  may intersect each other only at their ends. Moreover, these ends always  intersect either other arc or $\partial \Omega$.
   These intersections are called the  vertices
of the network.

A network is said to evolve by curve shortening flow if the evolution set is a network for every fixed time and each of its constituent arcs  $\gamma_i$ satisfies the shortening flow equation:
$$\frac{d \gamma_i}{dt}\cdot \nu=k_{\gamma_i}.$$

In order to have a well defined equation it is necessary to impose other  conditions at the interior nodes.
The most natural condition in the context of Brakke's  work (\cite{Br}) is to impose that the interior vertices are all trivalent and that the meeting angles of the curves are all equal to $\frac{2\pi}{3}$. When a  network satisfies this condition during its evolution we say that the network is regular.  We would like to remark that 
it is possible to impose other  conditions for the angles at the interior nodes, but these will not  be consider in this paper. At the exterior nodes we will impose Dirichlet boundary condition, namely we will prescribe the evolution of these nodes (mostly considering them to be fixed).

Recently, in \cite{sssnf} was proved that for an initial condition of $k$ half-lines there exists a connected  solution to the network flow.
Moreover, this solution is regular for all positive times and  self-similar; 
 however, it is not necessarily unique. In this paper we
 discuss uniqueness for these connected networks within a topological class.
 
To define the topological class we assume that the exterior
 nodes are fixed. If such vertices do not exist, that is if the curves
 extend to infinity, by  ''fixing the vertices" we mean that at infinity
 the network is asymptotic to certain fixed curves (which, in the case
 of the self-similar solutions described above, agree with the $k$
 half-lines that were taken as initial condition).

\begin{definition}
We say that two networks belong to the same topological class if there is a homotopy between them relative to  the exterior vertices.
\end{definition}

\begin{figure}[ht]
\subfigure[Two networks of the same topological class]{
\includegraphics[scale=0.4]{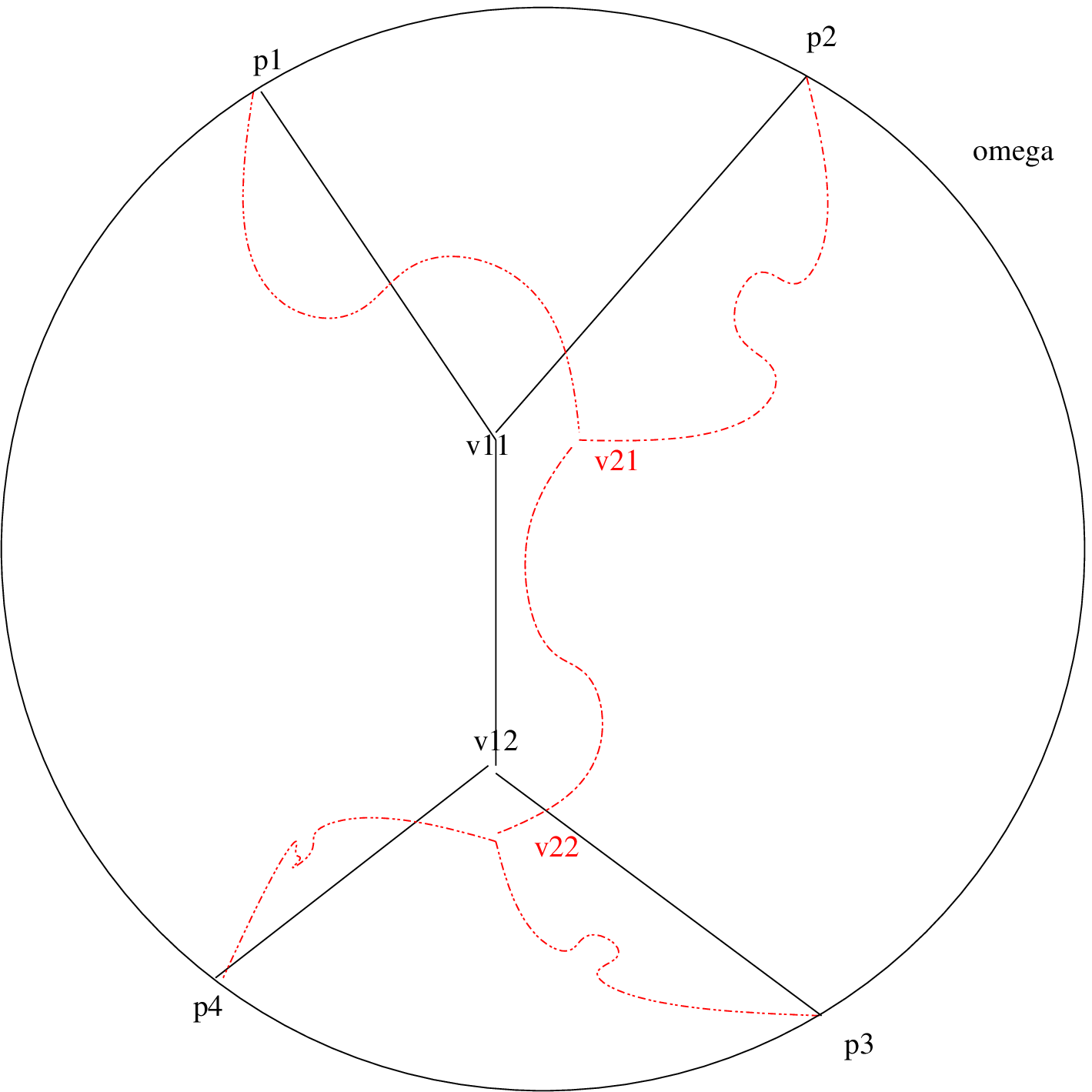}
\label{sametopclass}
}
\subfigure[Two networks of the different topological class]{
\includegraphics[scale=0.4]{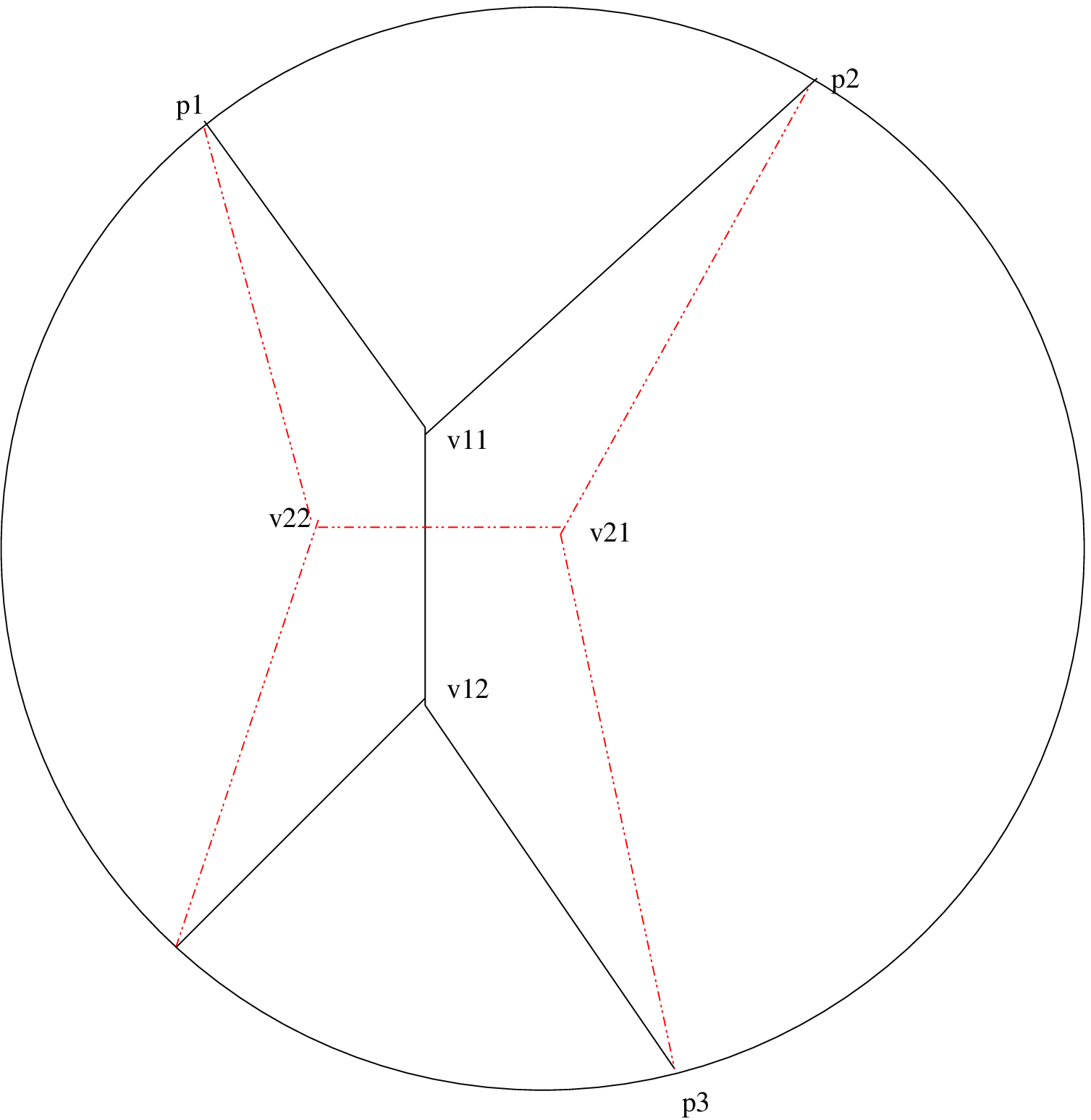}
\label{difftopclass}
}
\label{topclasses}
\caption[Topological Class]{Topological Class of networks with 4 exterior nodes}

\end{figure}

As an example of this concept in Figure \ref{sametopclass}  we show two networks in the same topological class, which contrast with Figure \ref{difftopclass}, where the two networks are in different topological classes. We remark that although in certain situations the networks in Figure  \ref{difftopclass} can be in the same topological class after rotation, we do not allow this by fixing the ends in the homotopy process.

\psfrag{p1}{$p_1$} 
\psfrag{p2}{$p_2$} 
\psfrag{p3}{$p_3$} 
\psfrag{p4}{$p_4$} 
\psfrag{v11}{$v_{11}$} 
\psfrag{v12}{$v_{12}$} 
\psfrag{v21}{$v_{21}$} 
\psfrag{v22}{$v_{22}$} 
\psfrag{omega}{$\Omega$} 

Now we can state the main theorem of this paper:

\begin{teo}\label{main theorem}
For each topological class there is at most one connected self-similar solution to the network flow with initial condition consisting of $k$ lines meeting at the origin and that is regular for positive times.
\end{teo}

This theorem can be extended as follows: \begin{teo}\label{main cor1}
Suppose that we have an arbitrary connected tree-like network that has
an evolution via network flow in a domain $\Omega$ that is regular for
$0<t<T$ and that at the boundary satisfy Dirichlet conditions.
Moreover, assume that the curvatures $k^i$ of this network
are bounded by $\frac{C}{\sqrt{t}}$.  
Then, given a topological
class, there is a unique evolution for this network.
\end{teo}

We would like to remark that for unbounded domains $\Omega$ we assume that the boundary conditions satisfied by the networks above are analogous to the ones satisfied by the self-similar solutions.
Namely, we assume that the arcs converge strongly (at least in $C^{2,1}$) to a fixed curve, which is assumed to be  compatible with the evolution equation (that in the case of the self-similar solutions described earlier, correspond to  lines and in general could be any unbounded curve that is a solution to the curve shortening flow equation).


The evolution of networks can be also seen from a different point of view, namely as the nodal set of the limit of solutions to the {\em vector-valued parabolic Allen-Cahn} equation. 
The Allen-Cahn Equation (with Dirichlet boundary condition) is given in a domain $\Omega$ by:
\begin{align}\frac{\partial u_\epsilon}{\partial t}-\Delta u_\epsilon+ \frac{\nabla_uW(u_\epsilon)}{\epsilon^2}&=0 \hbox{ for } x\in \Omega \label{laeq}\\
 u_\epsilon(x,0)&=\psi_\epsilon(x),\label{ci} \\
u_\epsilon|_{\partial \Omega}&=\phi_\epsilon (x,t) \label{bc} \end{align}
where $u_\epsilon:\rr^n\times \rr_+\to \rr^m$ and $W:\rr^n\to \rr$ is a positive potential with a finite number of minima. In particular we will concentrate on the case $m=n=2$ and $W$ is a function with 3 minima. 

In \cite{onthrlb} Bronsard and Reitich studied formally this equation
and predicted that as $\epsilon\to 0$, solutions would converge
a.e. to minima of $W$ and that the interfaces between these sets (where the solution converges to minima) might develop
a network structure that evolves under curve shortening flow. In
\cite{triodginz} was proved that this in fact holds for appropriate
potentials and initial conditions if the considered networks are {\em
triods} (networks that contain only three arcs). Moreover, it was
proved that any smooth evolution of triods can be realized as a nodal
set as described above. In order to prove Theorems \ref{main theorem}
and \ref{main cor1} we show that such a representation can be extended
to more general tree-like networks that have a regular evolution via
network flow for $0<t<T$.

We organize the paper as follows: In Section \ref{tcc} we show that the topological class of  network defines a unique coloring. In 
Section \ref{as} we show that this coloring gives us an Allen-Cahn approximation. In Section \ref{pmt} we use the previous approximation to conclude 
Theorems \ref{main theorem} and  \ref{main cor1}.

\begin{rem}
We would like to remark that along the coming proofs all constants will be denoted by $C$, but they might vary from line to line.
\end{rem}

\section{Topological class and coloring} \label{tcc}

\begin{figure}[ht]
 \epsfysize=3cm
\epsfig{ file=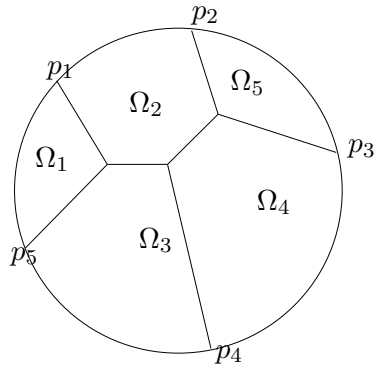,   width=0.3 \textwidth}
\caption{Partition of the domain induced by a network}
\label{partition}
\end{figure}

Notice that any network $\calN=\{\gamma_i\}$ contained in a set $\Omega$ defines a partition $\{\Omega_i\}$ of $\Omega\setminus\calN$ as follows: $\partial \Omega_i\subset
\calN$,  $\Omega_i$ is connected and $\bigcup_i\Omega_i=\Omega\setminus\calN$ (see as an example Figure \ref{partition}).

As stated in the introduction we are going to approximate solutions to the network flow via  solutions to the Allen-Cahn equation with a three well potential. 
Let us assume  that these three minima are given by $c_1,c_2$ and $c_3$. We will understand each of these minima as a color. It is expected that in each $\Omega_i$
the sequence of solutions $u_\epsilon$ to \equ{laeq} converge to one of the $c_i$. Hence,
a necessary requirement to approximate a network as an interface of the Allen-Cahn equation is that the partition defined by the network has a {\em three-coloring} as defined below:

\begin{definition}
Suppose that we have a domain $\Omega \in \rr^2$ and a partition $P=\{ \Omega_i \}_{i=1}^n$ of $\Omega$. Then we say that we can {\em three-color}
$P$ if for any region we can assign one of three fixed colors and the colors  assigned in any  adjacent regions  are different.
\end{definition}

An example of a three-coloring would be in Figure \ref{partition} to associate $c_1$ to $\Omega_1$ and $\Omega_4$, $c_2$ to $\Omega_2$ and $c_3$ to $\Omega_3$ and $\Omega_5$.

\psfrag{p1}{$p_1$} 
\psfrag{p2}{$p_2$} 
\psfrag{p3}{$p_3$} 
\psfrag{p4}{$p_4$} 
\psfrag{p5}{$p_5$} 
\psfrag{omega1}{$\Omega_1$} 
 \psfrag{omega2}{$\Omega_2$} 
 \psfrag{omega3}{$\Omega_3$} 
 \psfrag{omega4}{$\Omega_4$} 
 \psfrag{omega5}{$\Omega_5$} 


Suppose now that we have  a tree-like graph that all its interior nodes are trivalent and the exterior ones are simple. Then the following proposition holds:


\begin{prop} \label{coloring}
Suppose that we have a domain and a connected tree-like network with
only trivalent nodes in the interior and simple nodes in the
boundary. Let $P=\{\Omega_i\}$ be the partition associated to this
network. There is a unique three-coloring of $P$ (up to re-labeling
of the colors).  Moreover, the same coloring can be associated to any
other tree-like network that belongs to the same topological class.
\end{prop}

\begin{proof}
For simplicity we may assume is included in a compact domain like in Figure \ref{coloringfig}  (otherwise the exterior vertices that we will refer to correspond to ''vertices at infinity").

It is easy to see from the coming proof that the same coloring can be
associated to any two networks that belong to the same topological
class, hence we will not refer further to this assertion. For the
existence and uniqueness we will use induction on the number of
exterior vertices.

For $n=3$ the result is trivial (each region has one of the color and
this is unique, up to re-labeling of the colors). Suppose that the
result holds for networks of trivalent interior nodes and $n-1$ simple
exterior vertices. We prove that it also holds for networks of
trivalent interior nodes and $n$ simple exterior vertices. We will try
to illustrate the induction procedure in Figures \ref{coloring1}, \ref{coloring2} and \ref{coloring3} .  Let
us label the exterior vertices by $p_1,\ldots,p_n$.

We first  claim that there is an interior node $v$ and two exterior nodes $p_i$ and $p_j$ such that there are edges $e_{p_iv}$ and $e_{p_j v}$ that connect $p_i$ and $v$ and $p_j$ and $v$ respectively. This can be easily proved  by a combinatorial argument: Since the graph is a connected tree, we have that the number of interior nodes is equal to $n-2$ (this  computation can be found for example in \cite{sssnf}). 
Hence,  it is not possible that each exterior node is connected to a separate interior node.  Without loss of generality we can assume that the nodes adjacent to $v$ are $p_1$ and $p_2$.  And, since $v$ is a trivalent node, there is only one  remaining adjacent edge to $v$ that we label by $e$.

Now remove $p_1$, $v$ and the edge $e_{p_1v}$ and we regard $e_{p_2 v}\bigcup
e$ as one edge. Notice that by this procedure we eliminated only one region, which had $e_{p_1v}$ and $e_{p_2 v}$ as boundaries. 
Moreover,  we obtain a graph that topologically is equivalent to an tree-like graph with interior trivalent nodes and with $n-1$ single exterior nodes . By induction hypothesis we have that this graph can be uniquely colored with three colors (up to re-labeling of the colors). The edge that now is formed by $e_{p_2 v}\bigcup
e$ has two adjacent regions that are colored by 2 different colors, let us say $c_1$ and $c_2$. Now we can add again $p_1$, $v$ and the edge $e_{p_1v}$. This procedure divides one of the regions above in two and adds one region, that can be colored by $c_3$. Since this region is only neighboring the regions that for the previous graph were adjacent to $e_{p_2 v}\bigcup
e$, we have that this is a coloring of the graph with three colors.

\psfrag{p1}{$p_1$} 
\psfrag{p2}{$p_2$} 
\psfrag{p3...}{$p_3 \ldots$ } 
\psfrag{pn}{$p_n$} 
\psfrag{v}{$v$} 
\psfrag{e1v}{$e_{p_1v}$} 
\psfrag{e2v}{$e_{p_2v}$} 
\psfrag{e}{$e$} 
\psfrag{c1}{$c_1$} 
\psfrag{c2}{$c_2$} 
\psfrag{c3}{$c_3$} 
\begin{figure}[ht]
\centering
\subfigure[The original tree-like network]{
\includegraphics[scale=0.3]{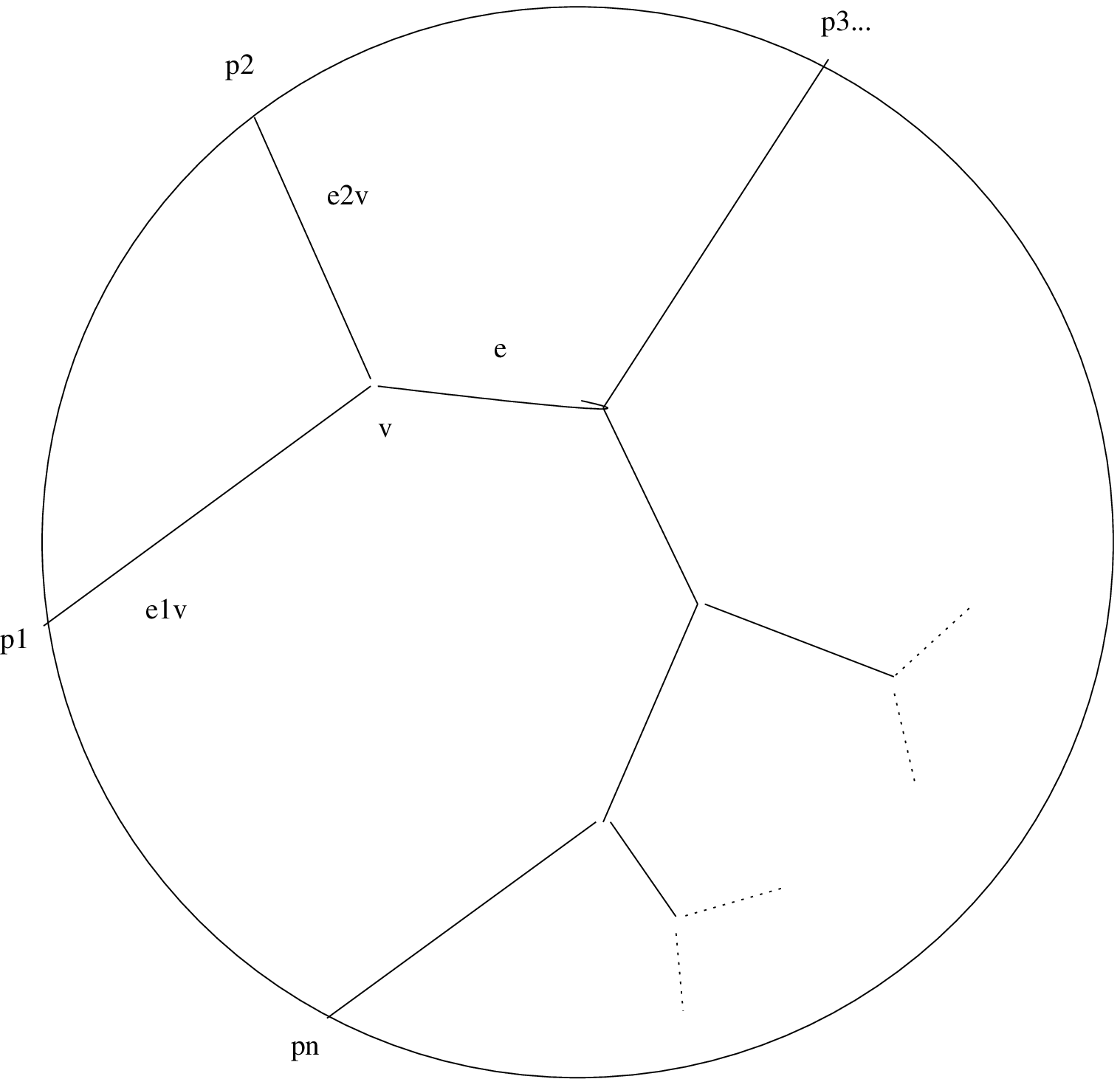}
\label{coloring1}
}
\subfigure[Removing  $e_{p_1v} $ and using the induction hypothesis]{
\includegraphics[scale=0.3]{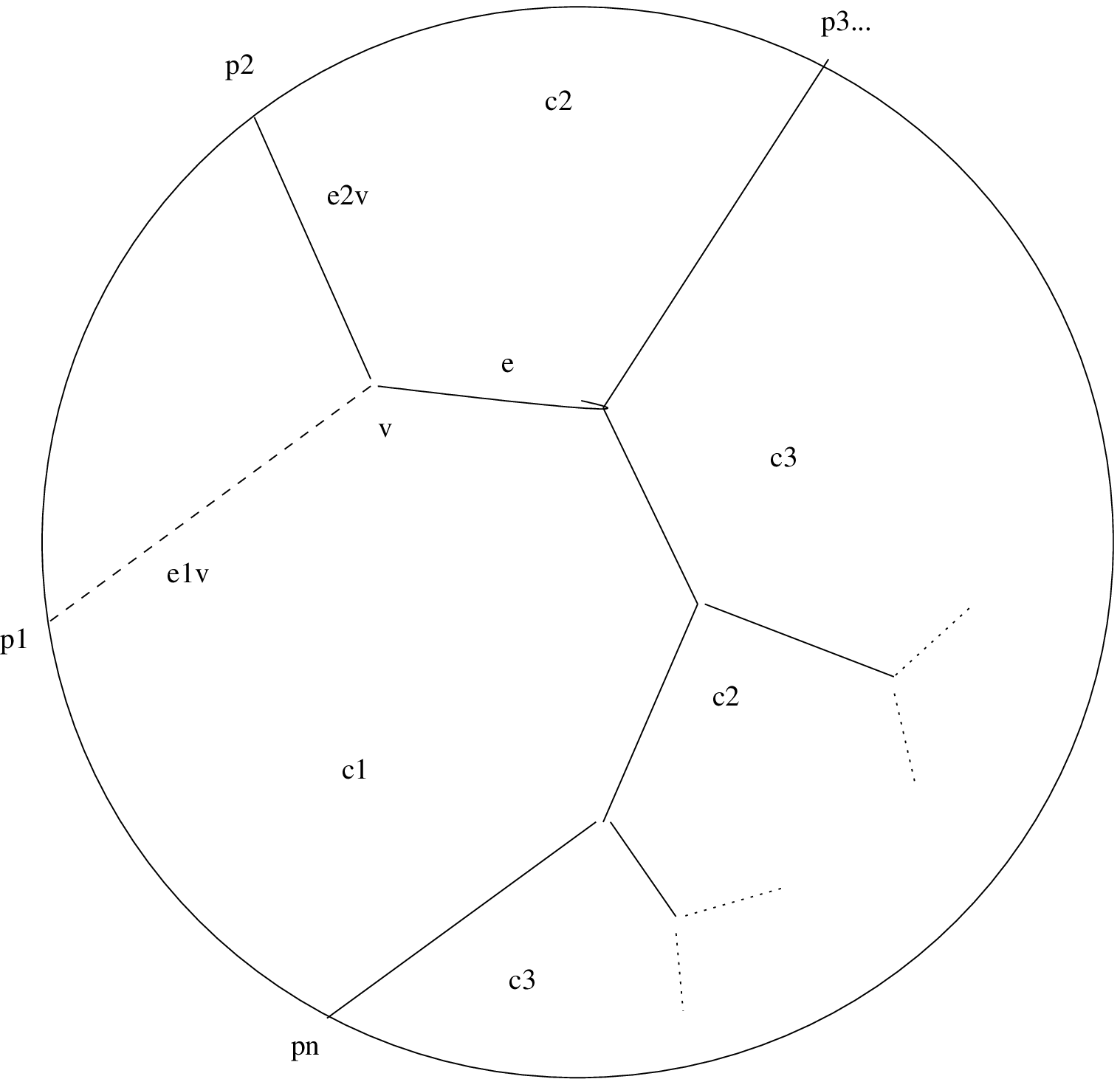}
\label{coloring2}
}
\subfigure[adding  $e_{p_1v}$ and coloring $R$]{
\includegraphics[scale=0.3]{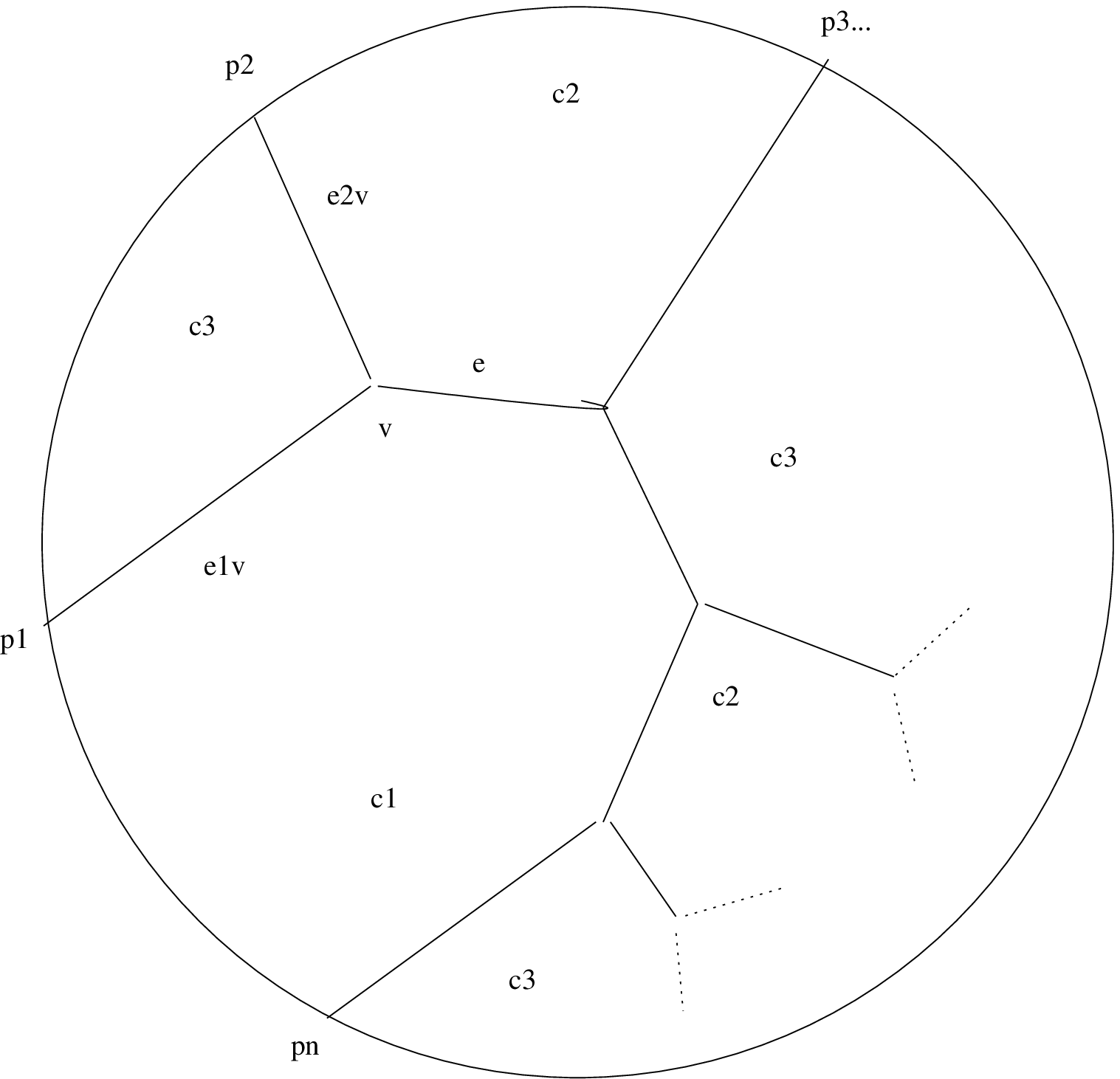}
\label{coloring3}
}
\label{coloringfig}
\caption[Proof of existence]{Proof by induction of the existence of a coloring}
\end{figure}

To prove the uniqueness we just invert the procedure. Suppose that we
have an arbitrary coloring of a network with $n$ exterior nodes. We
will show that this coloring agrees (up to re-labeling) with the one
just constructed above.  We can assume that $p_1, \ p_2,\ v,\ e,\
e_{p_1v}$ and $e_{p_2 v}$ are defined as before.  Via re-labeling we
can color the region enclosed by $e_{p_1v}$ and $e_{p_2 v}$ (that we
label by $R$) with $c_3$.  By removing $e_{p_1v}$, we have that this
region becomes part of a region of either color $c_1$ or $c_2$. Let us
assume it is $c_1$ and re-assign this color to region $R$ (and the
removed edge $e_{p_1v}$). Now, as before, we have a three-coloring of a
tree-like graph with interior trivalent nodes and with $n-1$ single
exterior nodes. By the uniqueness given by the induction hypothesis we
conclude that (re-labeling, if necessary) the coloring for the network
with the removed edge is the one used in the proof of existence
above. Moreover, since the regions that were adjacent to $R$ had color
$c_1$ and $c_2$ the only possibly necessary rearrengement would be to interchange 
$c_1$ and
$c_2$.  Hence, following the existence proof, we conclude that the
chosen arbitrary coloring (up to re-labeling) has to be the one
constructed previously.

\end{proof}

\section{Approximation skim}\label{as}

In this section we prove that in fact a regular tree-like network can
be understood as the nodal set of the limit of solutions to the
Allen-Cahn.  We assume that for positive times the nodes of the curves
are trivalent and that the meeting angles of the curves at the nodes
are all equal to $\frac{2\pi}{3}$. Although for certain networks it
would be possible to extend the proof for other fixed meeting angle, in
order to keep the presentation simpler, we will not discuss that
situation. In our case it is enough to consider {\em
symmetric} potentials 
$W$
with three minima (where $W$ attains the value 0). This corresponds to 
the following assumption
$\Gamma(c_i,c_j)=\Gamma(c_j,c_k)$ for every $i\ne j$, $j\ne k$, where
\begin{align*} \Gamma (\zeta_1,\zeta_2)=\inf\left\{\int_0^1W^{\frac{1}{2}}(\gamma(\lambda))|\gamma'(\lambda)|d\lambda: \right.&
 \gamma \in C^1([0,1],\rr^2), \notag \\
&\left. \gamma(0)=\zeta_1\hbox{ and }\frac{}{} \gamma(1)=\zeta_2 \right\}. \label{defdist}\end{align*}

Moreover, we assume that there are three unique heteroclinics
$\zeta_{ij}$ associated to this potential:{\em i.e} for each pair of
minima $c_i$ and $c_j$ of $W$ there is a unique curve $\zeta_{ij}$
that satisfies
\be \zeta_{ij}''(\lambda)+\frac{\nabla W(\zeta_{ij}(\lambda))}{2}=0 \label{eczeta},\ee
\be \lim_{\tau\to -\infty}\zeta_{ij}(\tau)=c_i, \quad 
\lim_{\tau\to \infty}\zeta_{ij}(\tau)=c_j, \label{explim}\ee

In \cite{stationary} was proved that under this conditions, if $W$ is
regular enough, there is a solution $u_*(x):\rr^2\to \rr^2$ to
$$-\Delta u_*+\nabla W(u_*)=0,$$ which satisfies as $r\to \infty$
\be u_*(r \cos\theta,r\sin\theta) \to c_i \hbox{ for } \theta \in 
[\theta_{i-1}, \theta_i], \label{convu*} \ee
\be  u_*(r \cos\theta_i,r\sin\theta_i)\to \zeta_{ij}(0). \label{convu*2}\ee
where the angles $\theta_i$ are given by the potential $W$ in the following manner:
$$ \frac{\sin \alpha_1}{\Gamma(c_2,c_3)}= \frac{\sin \alpha_2}{\Gamma(c_1,c_3)}= \frac{\sin \alpha_3}{\Gamma(c_1,c_2)}, $$
with $\alpha_i=\theta_{i+1}-\theta_i$.

In \cite{triodginz} it was proved:

\begin{teo}\label{mainteo1stpaper}
Let $\calT=\{\gamma^i\}$ be a triod evolving under curve shortening flow in a convex bounded domain $\Omega$  and let $O(t)$ be the point where the curves $\gamma_i(\cdot,t)$ meet.  Assume that the angles at which the curves $\gamma^i$ meet at  the point $O(t)$ are fixed and that the lengths of the curves   $\gamma_i$ stay bounded away from 0.
Consider a proper non-degenerate potential $W$ with three global minima $c_1, c_2$ and $c_3$.  Moreover, assume that this potential $W$ is consistent with the fixed angles at $O(t)$. Then there  are functions $\phi_\epsilon$,  $\psi_\epsilon$ and $v_\epsilon$ such that there is a solution $u_\epsilon$ to
\equ{laeq}-\equ{ci}-\equ{bc} and
\be\lim_{ \epsilon\to 0}v_\epsilon(x,t)\in \{c_1,c_2, c_3\} \hbox{ a.e., }\ee  
\be\{(x,t): \lim_{ \epsilon\to 0}v_\epsilon(x,t)\not\in \{c_1,c_2, c_3\} \}=\calT \hbox{ and }\ee
\be \lim_{\epsilon \to 0}\sup_{\Omega\times[0,T]}|u_\epsilon-v_\epsilon|(x,t)=0. \ee
\end{teo}

There are two  main ingredients on the proof  of this theorem:
 the construction
 of the function $v_\epsilon$ 
and Lemma \ref{cotaprinc}(that is stated below). In what follows we review these two
elements of the proof.

 


We start by  reviewing the construction of $v_\epsilon$. Two regions can be distinguished: close to $O(t)$ and away from this point.
In order to define these two regions, we consider $\tilde{\delta}$ (that needs to be chosen small enough)  and we take 
 a ball of radius $\tilde{\delta}$ around $O(t)$. Away from $O(t)$, namely on the complement of this ball we construct the function $\phi_\epsilon$ as follows:

 
Let $d_i(x,t)=dist(x,\gamma^i(\cdot,t))$ the signed distance of  a point $x\in \rr^2$ to the curve $\gamma^i(\cdot,t)$  (where
the signs of the distance functions need to be chosen appropriately) and consider $\tau^i(t)=(\cos \theta^i(t),\sin \theta^i(t))$ to be the unit tangents at  $O(t)$ .
Define the  sets: $$D_{ii}(t)=\{x\in \Omega: d_i(x,t)\leq \delta\}$$
$$D_{ii+1}(t)=\left\{x\in \Omega: d_i(x,t)\geq \frac{\delta}{2}, \quad d_{i+1}(x,t)\leq -\frac{\delta}{2} \hbox{ and }d_{ii+1}\geq\frac{\delta}{2} \right\}\subset S_{ii+1}.$$
Let $\xi^{ext}_{ij}(x,t)$ be a partition of unity associated to these sets, namely these functions satisfy 
$0\leq \xi^{ext}_{ij}\leq 1$ (where $j\in \{i,i+1\}$), $supp \ \xi^{ext}_{ij}\subset D_{ij}$ (where $supp$ denotes the support) and for every $x\in\bigcup_{i,j}D_{ij}$ holds $\sum_{i,j\in\{i,i+1\}}\xi^{ext}_{ij}(x,t)=1$.
 Then for $x\in  \bigcup_{i,j}D_{ij}$ we define 
 \begin{align} \phi_\epsilon(x,t)=
\sum_{i=1}^3 &\left( \xi^{ext}_{ii}\left(x,t\right)
\zeta_{ii+1}\left(
\frac{d_i(x,t)}{\epsilon}\right) +\xi^{ext}_{ii+1}\left(x,t\right)
c_i  \right). \label{defbc}\end{align} 

and $v_\epsilon(x,t)=\phi_\epsilon(x,t)$ for points $x$ away from $O(t)$. Moreover,
we use $\phi_\epsilon(x,t)$ as boundary  condition in \equ{bc}.

Furthermore, we extend the function $\phi_\epsilon$ to the whole domain $\Omega$ as follows:
$$\phi^\eta_\epsilon(x,t)\equiv \left(1-\eta_1\left(\frac{r(x,t)}{2\epsilon}+1-\frac{\tilde{\delta}}{2\epsilon}\right)\right)\phi_\epsilon(x,t),$$
where  $r(x,t)=|x-O(t)|$ and 
$\eta_1:\rr\to \rr$ is a function such that $\eta_1(x)\equiv 1$   when $|x|\leq\frac{1}{2}$ 
and $\eta_1(x)\equiv 0$ for $|x|\geq 1$.
\medskip

For points close to $O(t)$ we take $\frac{1}{2}<\rho<1$ and define the regions $B_{\tilde{\delta}}(O(t))\setminus B_{\epsilon^\rho}(O(t))$ and
$ B_{\frac{\epsilon^\rho}{2}}(O(t))$.
Consider now $x\in B_{\tilde{\delta}}(O(t))\setminus B_{\epsilon^\rho}(O(t))$ and take 
$$\tilde{\phi}_\epsilon(x,t)=
\sum_{i=1}^3\left(\xi^{int}_{2i}(\theta-\theta(t))
\zeta_{ii+1}\left(
\frac{d_i(x,t)}{\epsilon}\right)
\frac{}{}+\xi^{int}_{2i-1}(\theta-\theta(t))c_i\right),$$ 
where  $\theta(t)$ is the angle formed by the tangent $\tau^1(t)$ with the $x$-axis and
 $\{\xi^{int}_i\}_{i=1}^6$ is a partition of unity associated to  the following family of intervals:
$$\calA_{2i}=\left(\theta_i-\theta_{int}, \theta_i+\theta_{int}\right)$$
$$\calA_{2i+1}=\left(\theta_i+\frac{\theta_{int}}{2},\theta_{i+1}-\frac{\theta_{int}}{2}\right),$$
with $\theta_{int}$  an angle chosen appropriately small and 
and $\theta_i$  the angles given by \equ{convu*}. Then define
 $$\tilde{\phi}^\eta_\epsilon(x,t) \equiv
 \left(1-\eta_2\left(\frac{x-O(t)}{\epsilon^\rho}\right)\right)
 \tilde{\phi}_\epsilon(x,t),$$ where $\eta_2:\rr^2\to [0,1]$ is a
 function that satisfies $\eta_2(x)\equiv 1$ when
 $|x|\leq\frac{1}{2}$.

On the other hand the function $v_\epsilon$ in the region that is  $\epsilon^\rho $ close from the triple point  is given by the stationary solution  $u_*$  to \equ{laeq} that satisfies  \equ{convu*}  and  \equ{convu*2}. More specifically, within the ball of radius $\tilde{\delta}$ we consider the function
\begin{align}\tilde{v}_\epsilon(x,t)=&\tilde{\phi}^\eta_\epsilon(x,t)+\eta_2\left(\frac{x-O(t)}{\epsilon^\rho}\right)u_*\left( \frac{R_{\theta(t)}(x-O(t))}{\epsilon}\right),\label{close to the node 0}\end{align}
where
$R_\theta$ represents the rotation matrix by an angle $\theta$.

Finally, we let
\begin{align} v_\epsilon(x,t)=&\phi^\eta_\epsilon(x,t)
+\eta_1\left(\frac{r(x,t)}{2\epsilon}+1-\frac{\tilde{\delta}}{2\epsilon}\right)
\tilde{v}_\epsilon(x,t).\label{defv}\end{align}

The initial condition to be chosen in \equ{ci} is
\be \psi_\epsilon(x)= v_\epsilon(x,0).\label{defci}\ee  


Before discussing the second main element of the proof, a few remarks are necessary. Let 
 $$F_ {\epsilon}(h,\psi_\epsilon)=-\int_0^t \int_{\Omega}\calH_{\Omega}(x,y,t-s)\left(\frac{\nabla_u W(h+\phi^\eta_{\epsilon})}{\epsilon^2}+P \phi^\eta_\epsilon\right)(y,s)dyds$$
\be+ \int_{\Omega}\calH_{\Omega}(x,y,t)(\psi_\epsilon (y)-\phi^\eta_\epsilon(y,0))dy,\label{deffunc}\ee
where $\calH_{\Omega}$ denotes the heat kernel in  $\Omega$.
 
Notice that fixed points $h_\epsilon$ of this functional are solutions to the equation
\begin{align}
\frac{\partial h_\epsilon}{\partial t}-\Delta h_\epsilon
+\frac{\nabla_u W(h_\epsilon+\phi^\eta_\epsilon)}{2\epsilon^2}&=-P \phi^\eta_\epsilon\quad \hbox{ in }\Omega\label{ecpar}\\
h_\epsilon(x,t)&=0 \hbox{ on }\partial \Omega\label{cbpar}\\
h_\epsilon(x,0)&=\psi(x)-\phi^\eta_\epsilon(x,0).\label{cipar}
\end{align}

In particular, defining $u_\epsilon(x,t)=h_\epsilon(x,t)+\phi^\eta(x,t)$ we have $u_\epsilon(x,t)$ satisfies \equ{laeq}-\equ{ci}-\equ{bc}.

Now we can state the second  main tool used to prove  Theorem \ref{mainteo1stpaper}:

\begin{lem}[Lemma 4.1 in \cite{stationary}]\label{cotaprinc}

Fix $K>0$. Consider the sequences of continuous functions $\psi_n,
w_n $ satisfying $\sup |\psi_n|,\ \sup|w_n|\leq K$. Let
 $\epsilon_n\to 0$ and $T_n>0$. 
Assume in addition that for every $0<\epsilon<1$ holds
$\sup_{x\in \Omega, t\in [0,T]}|h_\epsilon|(x,t)\leq K$.
Then for each $\psi_n,\epsilon_n$ the functional $F_{\epsilon_n}$ has a unique fixed point $h_{\epsilon_n}$ and holds either
\begin{enumerate}
\item \label{mejorcaso}$\lim_{n \to \infty}\sup_{\Omega\times [0,T_n]}|w_n-h_{\epsilon_n}|\to 0$,
or
\item \label{peorcaso} there is a constant $C$, independent of ${\epsilon}_n$ and $T_n$
such that 
$$\sup_{\Omega\times [0,T_n]}|w_n-h_{\epsilon_n}|\leq C \sup_{\Omega\times [0,T_n]}
|F_{\epsilon_n}(w_n,\psi_n)-w_n|.$$

\end{enumerate}
\end{lem}

 The proof of Theorem \ref{mainteo1stpaper}  follows from showing that  $ \sup_{\Omega\times [0,T_n]}
|F_{\epsilon_n}(v_\epsilon,\psi_\epsilon)-v_\epsilon|\to0$ as $\epsilon\to 0$ and the uniform bounds (independent of $\epsilon$) of solutions  to  \equ{laeq}-\equ{ci}-\equ{bc}  proved in Lemma 2.2 of \cite{stationary}.

\begin{rem} \label{usingrcutoof}
Suppose that there is an $R$ such that $B_R\subset \Omega$. Then 
it is also possible in the definition of the functional $F_ {\epsilon} $ (given by \equ{deffunc}) to consider instead of $\phi^\eta$ the function $\chi_R(x)\phi^\eta$  where $\chi_R:\bar{\Omega}\to [0,1]$ is a function that satisfies $\chi_R(x)=0$ for $x\in B_{\frac{R}{2}}$ and $\chi_R(x)=1$ for $x\in \partial \Omega$. 

In this case we would also need to substitute $\phi^\eta$  by  $\chi_R(x)\phi^\eta$ in \equ{ecpar}, \equ{cbpar} and \equ{cipar}.
The proof of Lemma \ref{cotaprinc}  nor the proof of Theorem \ref{mainteo1stpaper} are altered by this modification.
\end{rem}

Now we   check that the  proof of Theorem \ref{mainteo1stpaper} extends for any tree-like network. More specifically we show that

\begin{teo}\label{mainteo1st}
Let $\calN=\{\gamma^i\}_{i=1}^{2n-3}$ be a network  evolving under curve shortening flow in a convex bounded domain $\Omega$   at times $t\in[0,T)$ and let $\{O_i(t)\}_{i=1}^{n-2}$ be its interior nodes.  Assume that the evolution is regular for all $0<t\leq T$ and that the lengths of the curves   $\gamma_i$ stay bounded away from 0. Moreover, we assume that the nodes stay a fixed distance away from the boundary and from each other and that their linear and angular speeds remain bounded.  The curvature of the arcs is assumed to be bounded by  $\frac{C}{\sqrt{t}}$ where $C$ is a fixed constant.

Consider a proper non-degenerate symmetric  potential $W$ with three global minima $c_1, c_2$ and $c_3$.Then there  are functions $\phi_\epsilon$,  $\psi_\epsilon$ and $n_\epsilon$ such that there is a solution $u_\epsilon$ to
\equ{laeq}-\equ{ci}-\equ{bc} that satisfies
\be\lim_{ \epsilon\to 0}n_\epsilon(x,t)\in \{c_1,c_2, c_3\} \hbox{ a.e., }\ee  
\be\{(x,t): \lim_{ \epsilon\to 0}n_\epsilon(x,t)\not\in \{c_1,c_2, c_3\} \}=\calN \hbox{ and }\ee
\be \lim_{\epsilon \to 0}\sup_{\Omega\times[0,T]}|u_\epsilon-n_\epsilon|(x,t)=0. \ee
Furhermore, $\lim_{ \epsilon\to 0}n_\epsilon(x,t)$ defines a coloring of the network $\calN$.
\end{teo}

We will only briefly sketch the proof of this theorem. For further details on the missing computations we refer the reader to \cite{triodginz}.

\begin{proof}

The proof of this theorem is completely analogous to the proof of Theorem \ref{mainteo1stpaper}.  

We start by constructing the function $n_\epsilon$:
We consider a coloring of the partition associated to the network and define the function $col:\Omega \setminus \calN(\cdot,t)\times [0,T)\to \{c_i\}_{i=1}^3$ as
$col(x,t)=c_i$ if the color associated to $x$ at time $t$ is $c_i$. Notice that as long as the curves do no not disappear during the evolution, the topological class will not change, hence, the coloring is in fact preserved in time. 

As before, we denote
$d_i(x,t)=dist(x,\gamma^i(\cdot,t))$ the signed distance of  a point $x\in \rr^2$ to the curve $\gamma^i(\cdot,t)$, where
the signs of the distance functions need to be chosen appropriately.
 
Analogous to the sets $D_{ij}$ defined in the proof of Theorem \ref{mainteo1stpaper}, we consider:
$$\calD_{i}(t)=\{x\in \Omega: d_i(x,t)\leq \delta\}$$
$$\calC_{i}(t)=\left\{x\in \Omega: col(x,t)=c_i \hbox{ and }
|d_j(x,t)|\geq \frac{\delta}{2} \hbox{  for every } j  \right \}.$$
Define a partition of unity $\{\xi^{ext}_{\calD_j},\quad \xi^{ext}_{\calC_j}\},$ associated to these sets.
We denote $Z_{\calD_i}(\cdot)=\zeta_{jk}(\cdot)$ if $\gamma_i\in\partial \calC_j \bigcap\partial \calC_k$.

 Then we can define for $x\in\bigcup_{j}\calD_{j}\cup \bigcup_{j}\calC_j$  the function 
\begin{align}
 \phi_\epsilon(x,t)=
\sum_{i} & \xi^{ext}_{\calD_i}\left(x,t\right)
Z_{\calD_i}\left(
\frac{d_i(x,t)}{\epsilon}\right) + \sum_{i}\xi^{ext}_{\calC_i}\left(x,t\right) col(x,t)
 . \label{defbc}\end{align} 

Here we would like to remark that the appropriate choice of sign of $d_i$ corresponds to have for $x\in\calC_j(t)$ that 
$Z_{\calD_i}\left(
\frac{d_i(x,t)}{\epsilon}\right)\to c_j$  as $\epsilon\to 0$. 

This function will correspond to the function $n_\epsilon$ away from the nodes. We also use this function as boundary  condition in \equ{bc}.

Now we extend the function $\phi_\epsilon$ to the whole domain $\Omega$ using  appropriate cut-off functions. Namely we take
$$\phi^\eta_\epsilon(x,t)\equiv \left( \sum_{i=1}^{n-2}\left(1-\eta_1\left(\frac{r_i(x,t)}{2\epsilon}+1-\frac{\tilde{\delta}}{2\epsilon}\right)\right)\right)\phi_\epsilon(x,t),$$
where  $r_i(x,t)=|x-O_i(t)|$ and 
$\eta_1:\rr\to \rr$ is a function  that satisfies $\eta_1(x)\equiv 1$   when $|x|\leq\frac{1}{2}$ 
and $\eta_1(x)\equiv 0$ for $|x|\geq 1$.
\medskip


Now we extend the construction close to the nodes $O_i(t)$. 
Since we have a coloring, around each node $O_i(t) $ we can do a construction analogous to 
\equ{close to the node 0} and consistent with the coloring.
 Let   $\tau^j_i(t)=(\cos \theta^j_i(t),\sin \theta^j_i(t))$ to be the unit tangents at  $O_i(t)$. 
In the region 
$\epsilon^\rho $ far from the triple point we consider a partition of unity $\{\xi^{int}_i\}_{i=1}^3$ associated to  the family of intervals $\{\calA_j\}_{j=1}^6$, where 
$$\calA_{2j}=\left(\frac{2j\pi}{3}-\theta_{int}, \frac{2j\pi}{3}+\theta_{int}\right)$$
$$\calA_{2j+1}=\left(\frac{2j\pi}{3}+\frac{\theta_{int}}{2},\frac{2(j+1)\pi}{3}-\frac{\theta_{int}}{2}\right),$$
and $\theta_{int}$ is an angle chosen appropriately small. We define $Z_{\calA_{2i}}(\cdot)=\zeta_{jk}(\cdot)$ if $col(x,t)=c_j$ for $x\in \calA_{2i+1}$ and 
$col(x,t)=c_k$ for $x\in \calA_{2i-1}$. We label the adjacent arcs to $O_i(t)$ by $\{\gamma^j_i\}_{j=1}^3$ and the corresponding distance
$d^j_i(x,t)=dist(x,\gamma^j_i(\cdot,t))$. Moreover, we assume that $\gamma^j_i\bigcap  \calA_{2j}\ne \emptyset$ (while  $\gamma^j_i\bigcap \calA_{2k}=
\emptyset$ for every $k\ne 2j$).
 Then we take the function
$$\tilde{\phi}^i_\epsilon(x,t)=
\sum_{j=1}^3\left(\xi^{int}_{2j}(\theta-\theta^i(t))
Z_{\calA_{2i}}\left(\frac{d^j_i(x,t)}{\epsilon}\right)
\frac{}{}+\xi^{int}_{2i-1}(\theta-\theta(t))col(x,t)\right),$$ 
where  $\theta^i(t)$ is the angle formed by the tangent $\tau^i_1(t)$ with the $x$-axis.
As before we extend the function $\tilde{\phi}^i_\epsilon$ to the whole domain  by multiplying by an appropriate cut-off function:
 $$\tilde{\phi}^i_\epsilon(x,t) \equiv \left(1-\eta_2\left(\frac{x-O_i(t)}{\epsilon^\rho}\right)\right) \tilde{\phi}^i_\epsilon(x,t),$$
with 
$\frac{1}{2}<\rho<1$ as before and  $\eta_2:\rr^2\to \rr$  is a smooth function such that $\eta_2(x)\equiv 1$ when $|x|\leq\frac{1}{2}$.

The function $n_\epsilon$ in the region that is  $\epsilon^\rho $ close from the triple point  is given by the stationary solution $u_*$   to \equ{laeq} that satisfies  \equ{convu*}  and  \equ{convu*2}. More precisely we take
\begin{align}\tilde{n}^i_\epsilon(x,t)=&\tilde{\phi}^i_\epsilon(x,t)+\eta_2\left(\frac{x-O_i(t)}{\epsilon^\rho}\right)u_*\left( \frac{R_{\theta^i(t)}(x-O_i(t))}{\epsilon}\right),\label{close to the node}\end{align}
where
$R_\theta$ represents the rotation matrix by an angle $\theta$.

 Finally we define
\begin{align} n_\epsilon(x,t)=&\phi^\eta_\epsilon(x,t)
+\sum_i \eta_1\left(\frac{r_i(x,t)}{2\epsilon}+1-\frac{\tilde{\delta}}{2\epsilon}\right)
\tilde{n}^i_\epsilon(x,t).\label{defv}\end{align}

The initial condition \equ{ci} is given by:
\be \psi_\epsilon(x)= n_\epsilon(x,0).\label{defci}\ee

 It is easy  to see that the computations in \cite{triodginz} (used to show
that $|F_{\epsilon_n}(v_\epsilon,\psi_\epsilon)-v_\epsilon|\to 0$ uniformly as $\epsilon \to 0$) can be extended to the function $n_\epsilon$ constructed above. We will not review this computation here. However an analogous computation will be presented in the coming proof .

Now, as in \cite{triodginz},  we can use Lemma  \ref{cotaprinc} and the uniform bounds provided by  Lemma 2.2 in \cite{stationary} to conclude the result.
 \end{proof}

Moreover, we extend this result to the following situation
\begin{teo}\label{generalization}
Let $\calN=\{\gamma^i\}_{i=1}^{2n-3}$ be a network  evolving under curve shortening flow in a convex bounded domain $\Omega$  and let $\{O_i(t)\}_{i=1}^{n-2}$ be its interior nodes.  Assume that the angles at which the curves $\gamma^i$ meet at  the points $O_i(t)$ are of $\frac{2\pi}{3}$ and that there 
is a constant $C$ such that the lengths of the curves   $\gamma_i$ and the distances between the interior nodes are bounded below by $C\sqrt{t}$. Moreover, we assume that the nodes stay a fixed distance away form the boundary and that  their linear and angular speeds are uniformly   bounded.  The curvature of the arcs is assumed to be bounded by  $\frac{C}{\sqrt{t}}$ where $C$ is a fixed constant.
Consider a proper non-degenerate symmetric  potential $W$ with three global minima $c_1, c_2$ and $c_3$.Then there  are functions $\phi_\epsilon$,  $\psi_\epsilon$ and $n_\epsilon$ such that there is a solution $u_\epsilon$ to
\equ{laeq}-\equ{ci}-\equ{bc} and
\be\lim_{ \epsilon\to 0}n_\epsilon(x,t)\in \{c_1,c_2, c_3\} \hbox{ a.e., }\ee  
\be\{(x,t): \lim_{ \epsilon\to 0}n_\epsilon(x,t)\not\in \{c_1,c_2, c_3\} \}=\calN \hbox{ and }\ee
\be \lim_{\epsilon \to 0}\sup_{\Omega}|u_\epsilon-n_\epsilon|(\cdot,t)=0 \hbox{ for every }t\in[0,T]. \ee
Furthermore, $\lim_{ \epsilon\to 0}n_\epsilon(x,t)$ defines a coloring of the network $\calN$.
\end{teo}

Before proving this Theorem we will need a weaker version of Lemma \ref{cotaprinc}:

\begin{lem}\label{cotaprinc2} 

Consider the functional $F_\epsilon$ in Lemma \ref{cotaprinc}.
Fix $K>0$. Consider the sequences of continuous functions $\psi_n,
w_n $ satisfying $\sup |\psi_n|,\ \sup|w_n|\leq K$. Let
 $\epsilon_n\to 0$ and $0<T_n\leq T$. 
Assume in addition that for every $0<\epsilon<1$ holds
$\sup_{x\in \Omega, t\in [0,T]}|h_\epsilon|(x,t)\leq K$.
Then for each $\psi_n,\epsilon_n$ the functional $F_{\epsilon_n}$ has a unique fixed point $h_{\epsilon_n}$ and holds either
\begin{enumerate}
\item \label{mejorcaso1}$\lim_{n \to \infty}\sup_{\Omega\times [0,T_n]}t |w_n-h_{\epsilon_n}|\to 0$,
or
\item \label{peorcaso1} there is a constant $C$, independent of ${\epsilon}_n$ and $T_n$
such that 
$$\sup_{\Omega\times [0,T_n]}t |w_n-h_{\epsilon_n}|\leq C \sup_{\Omega\times [0,T_n]}t
|F_{\epsilon_n}(w_n,\psi_n)-w_n|.$$

\end{enumerate}
\end{lem}

\begin{proof}
Let us assume that  neither  \equ{mejorcaso1} nor  \equ{peorcaso1}  hold. 
Then for every $n$ there  are $T_n,  \epsilon_n$ and $w_n$ such that 
$$\sup_{\Omega\times [0,T_n]}t |w_n-h_{\epsilon_n}|\geq n \sup_{\Omega\times [0,T_n]}t
|F_{\epsilon_n}(w_n,\psi_n)-w_n|.$$

Fix $\delta>0$. Since $w_n$ and $h_{\epsilon_n}$ are uniformly bounded, we have 
$\sup_{\Omega\times [\delta,T_n]}
|F_{\epsilon_n}(w_n,\psi_n)-w_n|
\frac{1}{\delta} \sup_{\Omega\times [\delta,T_n]} t
|F_{\epsilon_n}(w_n,\psi_n)-w_n|\to 0$. From Lemma \ref{cotaprinc} we have that 
$$\lim_{n\to\infty}\sup_{\Omega\times [\delta,T_n]} |w_n-h_{\epsilon_n}|\to 0.$$

Hence
\begin{align*}\lim_{n\to \infty}\sup_{\Omega\times [0,T_n]}t |w_n-h_{\epsilon_n}|\leq & \sup\{ T \sup_{\Omega\times [\delta,T_n]} |w_n-h_{\epsilon_n}|,\sup_{\Omega\times [0,\delta]}t |w_n-h_{\epsilon_n}|\}\\
 \leq & \sup_{\Omega\times [0,\delta]}t |w_n-h_{\epsilon_n}|\\ 
 \leq & C\delta.\end{align*}
Since $\delta$ is arbitrary we conclude that $\lim_{n\to \infty}\sup_{\Omega\times [0,T_n]}t |w_n-h_{\epsilon_n}|=0$, which contradicts that \equ{mejorcaso} does not hold.


\end{proof}

Now we can prove Theorem \ref{generalization}:

\begin{proof}[ Proof of Theorem \ref{generalization}]
The proof of this Theorem is analogous to the one of the previous one. However,
  some modifications both in the construction of the function $\mve$ and in the computation of the bounds that we use to conclude are necessary.

We would first like to remark that the hypotheses of the theorem allow
that some of the vertices at time $t=0$ split into several for
positive times. However, in the construction below we would like  to have
 the same number of nodes for all times. Hence, the considered
nodes $\{O_i(t)\}_{i=1}^{n-2}$ will be the ones that are different for
positive times and some of them might agree for t=0. 
For example, for the self-similar solutions coming out of $n$
half-lines that start at the origin, we have $n-2$ nodes
$\{O_i(t)\}_{i=1}^{n-2}$ that agree at $t=0$.

We start by discussing the function $\mve$ that will be used. The
notation will be as in the previous proof. We first  modify the
sets $\calD_i$ and $\calC_i$ as follows:

$$\calD^\epsilon_{i}(t)=\{x\in \Omega: d_i(x,t)\leq \frac{\epsilon^\rho}{2}\}$$
$$\calC^\epsilon_{i}(t)=\left\{x\in \Omega: col(x,t)=c_i \hbox{ and }
|d_j(x,t)|\geq \frac{\epsilon^\rho}{4}\hbox{  for every } j  \right \}.$$

Now  $\{\xi^{ext}_{\calD^\epsilon_j},\quad  \xi^{ext}_{\calC^\epsilon_j}\}$ denotes a partition of unity associated to these sets.
As before, we define $Z_{\calD^\epsilon_{i}}(\cdot)=\zeta_{jk}(\cdot)$ if $\gamma_i\in\partial \calC^\epsilon_j \bigcap \partial\calC^\epsilon_k$.

 Then we  define  the function $\phi_\epsilon $ for $x\in\bigcup_{j}\calD_{j}^\epsilon\cup \bigcup_{j}\calC_j^\epsilon$ (that will be used as the new boundary condition  in  \equ{bc})  by
\begin{align}
 \phi_\epsilon(x,t)=
\sum_{i} & \xi^{ext}_{\calD_i^\epsilon}\left(x,t\right)
Z_{\calD_i}^\epsilon\left(
\frac{d_i(x,t)}{\epsilon}\right) + \sum_{i}\xi^{ext}_{\calC_i^\epsilon}\left(x,t\right) col(x,t)
 . \label{defbc3}\end{align} 

The signs for $d_i$ are chosen appropriately as before and
 the following extension for $\phi_\epsilon$ are considered:
$$\phi_\epsilon^\eta(x,t)\equiv \left( \sum_{i=1}^{n-2}\left(1-\eta_2\left(\frac{x-O_i(t)}{\epsilon^\rho}\right)\right)\right)\phi_\epsilon(x,t),$$
where  
$\eta_2:\rr\to \rr$ is the cut-off function described in previous proof.
Now we define
\begin{align} n_\epsilon(x,t)=&\phi^\eta_\epsilon(x,t)
+\sum_{i=1}^{n-2}\eta_2\left(\frac{x-O_i(t)}{\epsilon^\rho}\right) u_*\left( \frac{R_{\theta^i(t)}(x-O_i(t))}{\epsilon}\right)
\label{defvft}\end{align}
and
the initial condition \equ{ci} is given by:
\be \psi_\epsilon(x)= n_\epsilon(x,0).\label{defcift}\ee
 



 Let
 $w_\epsilon(x,t)=\mve(x,t)-\phi_\epsilon^\eta(x,t)$. Correspondingly,
 if we define the functional as in Remark \ref{usingrcutoof}, we would
 take $w_\epsilon(x,t)=\mve(x,t)-\xi_R(x)\phi_\epsilon^\eta(x,t)$. If
 the latter choice were made, in all the computations below the
 function $\phi^\eta_\epsilon$ would need to be replaced by
 $\xi_R(x)\phi_\epsilon^\eta(x,t)$ and the same conclusions would
 hold. We leave this case to be checked by the reader.
 
Since $w_\epsilon(x,t)=0$ for $x\in\partial \Omega$, we can write
$$w_\epsilon(x,t)=\int_0^t \int_{\Omega}\calH_{\Omega}(x,y,t-s)
P w_\epsilon(y,s) dyds
+ \int_{\Omega}\calH_{\Omega}(x,y,t)w_\epsilon (y,0)dy.$$ Which implies  
\begin{align}(F_\epsilon(w_\epsilon,\psi_\epsilon)-w_\epsilon)(x,t)=&\int_0^t \int_{\Omega}\calH_{\Omega}(x,y,t-s)\left(-\frac{\nabla_u W(w_\epsilon+\phi^\eta_{\epsilon})}{\epsilon^2}(y,s)-P (\phi^\eta_\epsilon+w_\epsilon)(y,s)\right)dyds\notag \\
&+ \int_{\Omega}\calH_{\Omega}(x,y,t)(\psi_\epsilon (y)-\phi^\eta_\epsilon(y)-w_\epsilon(y,0))dy\notag \\
=&\int_0^t \int_{\Omega}\calH_{\Omega}(x,y,t-s)\left(-\frac{\nabla_u W(\mve)}{\epsilon^2}-P \mve\right)(y,s)dyds.
\label{eqvint}\end{align}

Hence, it is enough to prove that $$\sup_{\Omega\times[0,T)} t \left|\int_0^t \int_{\Omega}\calH_{\Omega}(x,y,t-s)\left(\frac{\nabla_u W(\mve)}{\epsilon^2}+P \mve \right)(y,s)dyds\right|\to 0\hbox{ as }\epsilon\to 0.$$

We start by briefly reviewing the computation away from the
nodes. This computation, for all times, is analogous to the ones
used to prove Theorem \ref{mainteo1stpaper} in \cite{triodginz}.
We bound: $$I_\epsilon^{ext}(x,t)=t \left|\int_0^t \int_{\Omega\setminus
\bigcup_i
B_{\epsilon^\rho}(O_i(s))}\calH_{\Omega}(x,y,t-s)\left(\frac{\nabla_u
W(\mve)}{\epsilon^2}+P \mve \right)(y,s)dyds\right|$$ for any
$x\in\Omega$ and $0\leq t\leq T$. This integral can be subdivided as
$I^{ext}_\epsilon\leq  I^{ext}_{\calD^\epsilon_i\setminus
(\calD^\epsilon_i\bigcap
\calC^\epsilon_i)}+I^{ext}_{\calC^\epsilon_i\setminus
(\calD^\epsilon_i\bigcap
\calC^\epsilon_i)}+I^{ext}_{\calD^\epsilon_i\bigcap
\calC^\epsilon_i}$, where
\begin{align*}A(x,y, t, s)=&\calH_{\Omega}(x,y,t-s)\left(\frac{\nabla_u W(\mve)}{\epsilon^2}+P \mve \right)(y,s),\\
 I^{ext}_{\calD^\epsilon_i\setminus (\calD^\epsilon_i\bigcap \calC^\epsilon_i)}=&t \left|\int_0^t \int_{\calD^\epsilon_i\setminus (\calD^\epsilon_i\bigcap \calC^\epsilon_i)}A(x,y,t,s)dyds\right| 
\\ I^{ext}_{\calC^\epsilon_i\setminus (\calD^\epsilon_i\bigcap \calC^\epsilon_i)}=&y\left|\int_0^t \int_{\calC^\epsilon_i\setminus (\calD^\epsilon_i\bigcap \calC^\epsilon_i)}A(x,y,t,s)dyds\right| 
 \\I^{ext}_{\calD^\epsilon_i\bigcap \calC^\epsilon_i}=&t\left|\int_0^t \int_{\calD^\epsilon_i\bigcap \calC^\epsilon_i}A(x,y,t,s)dyds\right| 
\end{align*}

Using the definition of $\mve$ it is not hard to compute that
$$\begin{array}{rll}
 |A(x,y, t, s)|=&\left|\calH_{\Omega}(x,y,t-s) \frac{k_{i}^2(\lambda,t)d_{i}}{\epsilon(1+k_{i}(\lambda,t)d_{i})}\zeta'_{ii+1}\right| & \\
 \leq &
\calH_{\Omega}(x,y,t-s)\left| \frac{k_{i}^2(\lambda,t)d_{i}}{\epsilon(1+k_{i}(\lambda,t)d_{i})}\right| e^{-c\frac{d_i(x,t)}{\epsilon}}  & \hbox{ for } x\in \calD^\epsilon_i\setminus (\calD^\epsilon_i\bigcap \calC^\epsilon_i).
\\  & & \\
  |A(x,y, t, s)|=&0 &\hbox{ for } x\in \calC^\epsilon_i\setminus
  (\calD^\epsilon_i\bigcap \calC^\epsilon_i). \\ & & \\ |A(x,y, t,
  s)|=&\left|\calH_{\Omega}(x,y,t-s) \right|
  \left|P\mve+\frac{\nabla_u W\left( \mve\right)}{\epsilon^2}\right|&
  \\  &\leq C \max \left\{\sup_i | k_i|,
  \frac{e^{-c\epsilon^{\rho-1}}}{\epsilon^2}\right\} & \hbox{ for }
 x\in \calD^\epsilon_i\bigcap \calC^\epsilon_i.\\ & & \\
  \left|P\mve+\frac{\nabla_u W\left( \mve\right)}{\epsilon^2}\right| &
  \to 0 \hbox{ as } \epsilon \to 0 &\hbox{ for } x\in
  \calD^\epsilon_i\bigcup \calC^\epsilon_i.
\end{array}$$

Since $k^i$ is an integrable function, the Dominated Convergence Theorem implies that 
\be \lim_{\epsilon \to 0}I^{ext}_\epsilon=0.\label{compdeps}\ee
For further details in this computation we refer the reader to \cite{triodginz}.

\medskip

 Now we need to bound
 $$I_\epsilon^{int}(x,t)=t \left|\int_0^t \int_{\bigcup_i B_{\epsilon^\rho}(O_i(s))}\calH_{\Omega}(x,y,t-s)\left(\frac{\nabla_u W(\mve)}{\epsilon^2}+P \mve \right)(y,s)dyds\right|.$$
 
Notice that there is a $K$ that depends only on the network such that   for $t\geq K\epsilon^{2\rho}$ all the nodes are already separated enough to make the definition of $\mve$ appropriate to distinguish each region. It is not hard to check, that the computations in \cite{triodginz} can be easily extended for these large times. 
However, since at $t=0$ there might be nodes that coincide (such as in the case of the self-similar solutions described in the introduction), for  $t\leq K\epsilon^{2\rho}$ the computations need to be handled more carefully near this multiple nodes.
In order
  to have some room in our calculations, we will divide the times into 
slightly larger regions, namely $t\leq \epsilon^\rho$ and  $t\geq \epsilon^\rho$.
We will first focus on bounding $I_\epsilon^{int}(x,t)$ for small times.

 Since the nodes $O_i(t)$ are away form the boundary, we have $$
 \calH_{\Omega}(x,y,t-s)\left|\frac{\nabla_u W(\mve)}{\epsilon^2}+P
 \mve \right|\leq \frac{C}{\epsilon^2}
 \frac{e^{-\frac{|x-y|^2}{t-s}}}{t-s} .$$ 
 Then it follows for $t\leq \epsilon^\rho$:

 \begin{align*}I_\epsilon^{int}(x,t) \leq &\frac{C}{\epsilon^2}
t \int_0^{t} \int_{\bigcup_i B_{\epsilon^\rho}(O_i(s))}
 \frac{e^{-\frac{|x-y|^2}{t-s}}}{t-s} dyds \\ 
\leq& \frac{C \epsilon^\rho}{\epsilon^2}
 \int_0^{t-\epsilon^3}\int_{B_{\epsilon^\rho}}
 \frac{1}{t-s}   dyds,  + \int_{t-\epsilon^3}^t\int_{\rr^2}
 \frac{e^{-\frac{|x-y|^2}{t-s}}}{t-s}  dyds \\ \leq &
  \frac{C \epsilon^\rho}{\epsilon^2}\left( \pi \epsilon^{2\rho}(\ln t-3\ln\epsilon) + \epsilon^3\right)
 \\ \leq & -C\epsilon^{3\rho-2}\ln \epsilon
\end{align*}

Hence we have,
 \be I_\epsilon^{int}(x,t) \leq -C\epsilon^{3\rho-2}\ln \epsilon, \hbox{ for }
 t\in [0,  \epsilon^\rho].\label{short time at nodes}\ee
Consequently,   we choose $\frac{1}{2}< \frac{2}{3}<\rho<1$.
 
For $t\geq  \epsilon^{\rho}$ we bound the integral as follows:
$$I_\epsilon^{int}\leq T( I_{s\leq K \epsilon^{2\rho}}^{int} +  I_{s\geq K\epsilon^{2\rho}}^{int}+ I_{s\geq K\epsilon^{2\rho}}^{trans}), $$
where
\begin{align*}   
I_{s\leq K\epsilon^{2\rho}}^{int} =&\left|\int_0^{K\epsilon^{2\rho}}\int_{\bigcup_i B_{\epsilon^\rho}(O_i(s))}\calH_{\Omega}(x,y,t-s)\left(\frac{\nabla_u W(\mve)}{\epsilon^2}+P \mve \right)(y,s)dyds\right|\\
I_{s\geq K\epsilon^{2\rho}}^{int} =&\left|\int_{K\epsilon^{2\rho}}^t \int_{\bigcup_i B_{\frac{\epsilon^\rho}{2}}(O_i(s))}\calH_{\Omega}(x,y,t-s)\left(\frac{\nabla_u W(\mve)}{\epsilon^2}+P \mve \right)(y,s)dyds\right| \\
I_{s\geq K\epsilon^{2\rho}}^{trans} =&\left|\int_{K\epsilon^{2\rho}}^t \int_{\bigcup_i B_{\epsilon^\rho}(O_i(s))\setminus B_{\frac{\epsilon^\rho}{2}}(O_i(s))}\calH_{\Omega}(x,y,t-s)\left(\frac{\nabla_u W(\mve)}{\epsilon^2}+P \mve \right)(y,s)dyds\right|
.\end{align*}

For $s\leq K\epsilon^{2\rho}$,  as in the computation for \equ{short time at nodes}, we have 
$$ \calH_{\Omega}(x,y,t-s)\left|\frac{\nabla_u W(\mve)}{\epsilon^2}+P \mve \right|(y,s)\leq  \frac{C}{\epsilon^2} \frac{e^{-\frac{|x-y|^2}{t-s}}}{t-s} .$$
Then
\begin{align*} I_{s\leq K\epsilon^{2\rho}}^{int}(x,t) \leq &\frac{C}{\epsilon^2} \int_0^{K\epsilon^{2\rho}}\int_{\bigcup_i B_{\epsilon^\rho}(O_i(s))} \frac{e^{-\frac{|x-y|^2}{t-s}}}{t-s} dyds
\notag \\ \leq &\frac{C}{\epsilon^2}\left( \int_0^{K\epsilon^{2\rho}} \int_{\bigcup_i B_{\epsilon^\rho}(O_i(s))} \frac{1}{t-s} dyds \right)
 \notag \\ \leq &
-\frac{C}{\epsilon^2} \pi \epsilon^{2\rho} \sum_i\ln(1-\frac{K\epsilon^{2\rho}}{t})
\notag \\ \leq &-C \frac{\sum_i\ln(1-K\epsilon^{\rho})}{\epsilon^{2(1-\rho)}}
\notag \\ \leq &C(n-2) \epsilon^{3\rho-2}.
\end{align*}
Since we took  $\frac{1}{2}< \frac{2}{3}<\rho<1$,  we have
\be I_{s\leq K\epsilon^{2\rho}}^{int}(x,t)\leq C\epsilon^{\tilde{\rho}} \hbox{ for some }\tilde{\rho}>0  \label{1 bound at big t}\ee
 For $
  I_{s\geq K\epsilon^{2\rho}}^{int}$ we proceed as in \cite{triodginz}. It is easy to compute that
  
  \begin{align*}
|A(x,y, t, s)|=&
\left|\calH_{\Omega}(x,y,t-s)  Ju_*\frac{\theta' R_{\theta(t)+\frac{\pi}{2}}(x-O(t))-R_{\theta(t)}O'(t)}{\epsilon} \right| \\
\leq &
\frac{C}{\epsilon}\calH_{\Omega}(x,y,t-s) \hbox{ for } y \in B_{\frac{\epsilon^\rho}{2}}(O_i(t))
 \end{align*}
 
 where, as before, $A(x,y, t, s)=\calH_{\Omega}(x,y,t-s)\left(\frac{\nabla_u W(\mve)}{\epsilon^2}+P \mve \right)(y,s)$.
 
Since the nodes are away from the boundary, we have that
 $\left|\calH_{\Omega}(x,y,t-s)\right|\leq \frac{Ce^{-\frac{|x-y|^2}{t-s}}}{t-s}$ and
\begin{align}I_{s\geq K\epsilon^{2\rho}}^{int} 
\leq & \frac{C}{\epsilon}\int_{K\epsilon^{2\rho}}^t \int_{B_{\frac{\epsilon^\rho}{2}}(O(s))} \frac{e^{-\frac{|x-y|^2}{t-s}}}{t-s} dy ds\notag \\
\leq &\frac{C}{\epsilon}\left(\int_{K\epsilon^{2\rho}}^{t-\epsilon^m} \int_{B_{\frac{\epsilon^\rho}{2}}(O(s))} \frac{1}{t-s} dy ds+\int_{t-\epsilon^m}^t \int_{\rr^2} \frac{e^{-\frac{|x-y|^2}{t-s}}}{t-s} dy ds \right)\notag \\
\leq &\frac{C}{\epsilon}\left(\int_{K\epsilon^{2\rho}}^{t-\epsilon^m} \frac{\pi \epsilon^{2\rho}}{4(t-s)}ds+ \int_{t-\epsilon^m}^t ds \right)\notag \\ \leq &\frac{C}{\epsilon}\left((\ln( t-K\epsilon^{2\rho})-m\ln\epsilon)\frac{\pi \epsilon^{2\rho}}{4}
+ \epsilon^m\right)\notag \\ \leq &C\left((\ln T-m\ln\epsilon)\frac{\pi \epsilon^{2\rho-1}}{4}
+ \epsilon^{m-1}\right)\to 0 \hbox{ as } \epsilon \to 0. \label{bd1}\end{align}

For $t\geq \epsilon ^{\rho}$  in $B_{\epsilon^\rho}(O_i(t)) \setminus B_{\frac{\epsilon^\rho}{2}}(O_i(t))$
 we have that $\eta_2\left(\frac{x-O_j(t)}{\epsilon^\rho}\right) =0$ and $\mve= \phi_\epsilon(x,t)
+\eta_2\left(\frac{x-O_i(t)}{\epsilon^\rho}\right) \left(u_*\left( \frac{R_{\theta^i(t)}(x-O_i(t))}{\epsilon}\right)-\phi_\epsilon(x,t)\right)$. Hence

\begin{align}Pv_\epsilon+\frac{\nabla_u W(v_\epsilon)}{\epsilon^2}=&\left(P \phi_\epsilon(x,t)
+\frac{\nabla_u W(\phi_\epsilon)}{\epsilon^2}\right)\notag \\& +P\left[\eta_2\left(\frac{x-O_i(t)}{\epsilon^\rho}\right)
\left( u_*\left(\frac{R_{\theta^i(t)}(x-O_i(t))}{\epsilon} \right)-\phi_\epsilon(x,t)\right)\right]
\notag \\ &+\frac{\nabla_u W(v_\epsilon)}{\epsilon^2}-\frac{\nabla_u W(\phi_\epsilon)}{\epsilon^2} .\label{ect1}   \end{align}
 
The first parenthesis can be bounded as \equ{compdeps}, while the bounds for the other quantities follow from the fast convergence of $u_*$ to $\phi_\epsilon$
(see \cite{triodginz} and \cite{stationary} for details on this). 
Hence we have
 for any $m>0$

\be I_{s\geq K\epsilon^{2\rho}}^{trans}\leq  C\left(\int_0^t \int_{ B_{\epsilon^\rho}(O(t)) \setminus B_{\frac{\epsilon^\rho}{2}(O(t))}}
  \calH_{\Omega}(x,y,t-s)\left|P\phi_\epsilon+\frac{\nabla_u W\left( \phi_\epsilon\right)}{\epsilon^2}\right|(y,s)+\epsilon^{m}\right) . \label{bdt2}\ee

We conclude  from Lemma \ref{cotaprinc2} and \equ{compdeps}, \equ{short time at nodes},  \equ{1 bound at big t}, \equ{bd1} and \equ{bdt2} that for every $t>0$
the result holds. However, since the initial condition was appropriately chosen, we in fact have that the result holds for every $t\geq 0$.

\end{proof}

\begin{rem}
Notice that the self-similar solutions to the network flow studied in
\cite{sssnf} satisfy the conditions of  Theorem \ref{generalization} in a large ball for
a fixed time (that depends on the size of the domain). The boundary
condition needs to be chosen as the motion of the intersection point
of the solution and the boundary of the ball.
\end{rem}

\section{Proof of Theorems  \ref{main theorem} and \ref{main cor1}
} \label{pmt}

   Theorems \ref{main theorem} is a particular case of Theorem \ref{main cor1} . 
 Hence, we only treat the more general case.
 Suppose that we have two networks 
$\calN_1$ and $\calN_2$ in a domain $\Omega$. Furthermore assume that they  have the same topological type and their initial conditions coincide.

Since Theorem \ref{generalization} was proved in a bounded domain,  if $\Omega$ is unbounded we define the bounded domain $\Omega_R$ as the intersection of 
a large ball of Radius $R$  with $\Omega$ and we restrict the evolution to this domain. The requirements on the size of the ball will be specified later in this proof and the time $T$, up to when the evolution is considered,  is such that the nodes of the network stay away from the boundary of $\Omega_R$. 
Notice that if $\Omega$ is bounded domain and $R$ is large enough the same definition can be made without altering $\Omega$. Hence we will always refer to the domain as $\Omega_R$.
 
Since the topological classes  of $\calN_1$ and $\calN_2$ agree, Proposition \ref{coloring} implies that the same coloring can be associated to both of them.  
From Theorem \ref{generalization} that there are solutions  $u_\epsilon^1$ and $u_\epsilon^2$ to \equ{laeq}-\equ{ci}-\equ{bc} which as $\epsilon\to 0$ develop nodal sets that respectively agree with $\calN_1$ and $\calN_2$ in $\Omega_R$ and respect the associated  coloring. 
Let $ \phi^1_\epsilon$ and $ \phi^2_\epsilon$ be the boundary conditions of $u_\epsilon^1$ and $u_\epsilon^2$ respectively. To simplify the notation, we assume that these functions have been already extended to the whole domain via a cut-off function. 
We denote
$$\tilde{\phi}_\epsilon^i(x,t)=\left\{ 
\begin{array}{ll}\chi_R(x)\phi^i_\epsilon(x,t) & \hbox{ if } \Omega \hbox{ is bounded}\\ \phi^i_\epsilon(x,t) & \hbox{ otherwise } \end{array}\right. .$$

 Furthermore Lemma \ref{cotaprinc} implies (by choosing
 $w_\epsilon(x,t)=u^2_{\epsilon}(x,t)-\tilde{\phi}_\epsilon^1(x,t)$)
 that either $$\sup_{\Omega_R\times
 [0,T_n]}|u^1_\epsilon-u^2_{\epsilon}|\to 0$$ or $$\sup_{\Omega_R\times
 [0,T_n]}|u^1_\epsilon-u^2_{\epsilon}|\leq C \sup_{\Omega_R\times
 [0,T_n]}
 |F_{\epsilon}(u^2_\epsilon-\tilde{\phi}^1)-u^2_\epsilon+\tilde{\phi}^1|,$$
 where
\begin{align*}F_ {\epsilon}(h)=& -\int_0^t \int_{\Omega_R}\calH_{\Omega_R}(x,y,t-s)\left(\frac{\nabla_u W(h+\tilde{\phi}^1_{\epsilon})}{\epsilon^2}+P ( \tilde{\phi}^1_\epsilon)\right)(y,s)dyds\\
&+ \int_{\Omega_R}\calH_{\Omega_R}(x,y,t)(\psi_\epsilon (y)-\tilde{\phi}^1_\epsilon(y,0))dy,\label{deffunc}\end{align*}



Let us assume that we are in the second case.
Since $u_2$ is a solution to \equ{laeq} and $ u^2_\epsilon(x,t)-\tilde{\phi}^2_\epsilon(x,t)=0$ for $x\in\partial \Omega_R$ we have that
\begin{align*} u^2_\epsilon(x,t)-\tilde{\phi}^2_\epsilon(x,t)= &-\int_0^t \int_{\Omega_R}\calH_{\Omega}(x,y,t-s)\left(\frac{\nabla_u W(u^2_\epsilon)}{\epsilon^2}+P (\tilde{\phi}^2_\epsilon)\right)(y,s)dyds\\
&+ \int_{\Omega_R}\calH_{\Omega}(x,y,t)(\psi_\epsilon (y)-\tilde{\phi}^2_\epsilon(y,0))dy,
 \end{align*} 
 and
 
 Notice that the functions $\psi_\epsilon$ coincide since the initial condition of the networks agree.
   Hence, we have:


 \begin{align*} ( F_{\epsilon}(u^2_\epsilon-
 \tilde{\phi}_\epsilon^1)-u^2_\epsilon+\tilde{\phi}_\epsilon^1)(x,t)=&
 -\int_0^t \int_{\Omega_R}\calH_{\Omega}(x,y,t-s)\left(P (
 \tilde{\phi}^1_\epsilon-\tilde{\phi}^2_\epsilon)\right)(y,s)dyds \\ &
 +
 \int_{\Omega_R}\calH_{\Omega}(x,y,t)(\tilde{\phi}^1_\epsilon(y,0)-\tilde{\phi}i^2_\epsilon(y,0))dy.
\end{align*}

Since, as $R\to \infty$ we have that $ |\tilde{\phi}^1_{\epsilon}-\tilde{\phi}^2_{\epsilon}|_{C^{2,1}}\to 0$, we have that for any $\delta>0$ there is an $R$ such that
\be\sup_{(x,t)\in\Omega_R\times[0,T]}\left| F_{\epsilon}(u^2_\epsilon- \tilde{\phi}_\epsilon^2)-u^2_\epsilon-\tilde{\phi}_\epsilon^2\right|<\delta. \label{finalbound1}\ee

Now, if the evolutions of $\calN_1$ and $\calN_2$ are different we would have that $\sup_{\Omega\times [0,T_n]}|u^1_\epsilon-u^2_{\epsilon}|>\min_{i\ne j}\{|c_i-c_j|\}$. Choose $\delta$ such that $C\delta< \min_{i\ne j}\{|c_i-c_j|\}$ (where $C$ is the constant given by Lemma \ref{cotaprinc}). Then Lemma \ref{cotaprinc} and \equ{finalbound1} imply:
$$\min_{i\ne j}\{|c_i-c_j|\}<\sup_{\Omega\times [0,T_n]}|u^1_\epsilon-u^2_{\epsilon}|\leq C \sup_{\Omega\times [0,T_n]}\leq \sup_{\Omega\times [0,T_n]}| F_{\epsilon}(u^2_\epsilon- \tilde{\phi}_\epsilon^2)-u^2_\epsilon+\tilde{\phi}_\epsilon^2|\leq C\delta,$$
which yields a contradiction.
$\Box$

\begin{rem}
We would like to remark that in the previous proof we assumed that the constant $C$ in Lemma \ref{cotaprinc} can be chosen uniformly for every domain $\Omega_R$. This in fact holds,  since the bounds for  $\tilde{\phi}_\epsilon$ and for $u_\epsilon$ provided by Lemma 2.2 in \cite{stationary} are uniform in $\rr^2$.
\end{rem}




\end{document}